\documentstyle[12pt]{amsart} 
\makeatletter
\def\displaylines#1{\displ@y
  \halign{\hbox to\displaywidth{$\@lign\hfil\displaystyle##\hfil$}\crcr
    #1\crcr}}
\makeatother

\numberwithin{equation}{section}
\newcommand{\beq}{\begin{equation}}
\newcommand{\eeq}{\end{equation}}

\newtheorem{theorem}{Theorem}[section]

\newtheorem{prop}[theorem]{Theorem}

\newtheorem{coroll}[theorem]{Corollary}

\newtheorem{lemm}[theorem]{Lemma}

\newtheorem{anoprop}[theorem]{Proposition}

\newtheorem{eg}{\rm\sl \uppercase{Example}}[section]

\begin{document}

\title{Relative Cohomology of Banach Algebras}
\author{Zinaida A. Lykova}
\address{Zinaida A. Lykova\\
 Department of Mathematics and Statistics\\
 Lancaster University \\
 Lancaster LA1 4YF \\
 England}
\email{maa032@@cent1.lancs.ac.uki}

\begin{abstract}
     Let $A$ be a Banach algebra, not necessarily unital, and let $B$ be a 
closed subalgebra of $A$. We establish a connection between the Banach 
cyclic cohomology group $ {\cal HC}^n(A)$  of $A$ and the Banach 
$B$-relative cyclic cohomology group $ {\cal HC}^n_B(A) $ 
 of $A$. We prove that, for a Banach algebra $A$ with a bounded approximate 
identity and an amenable closed subalgebra $B$ of $A$, up to topological
 isomorphism, ${\cal HC}^n(A) = {\cal HC}^n_B(A) $ for all $n \ge 0$. We 
also establish a  connection between the Banach simplicial or  
cyclic cohomology groups of $A$ and those of the quotient algebra $A/I$ by an 
amenable closed bi-ideal $I$. The results are applied to the calculation of these 
groups for certain operator algebras, including von Neumann algebras.

\end{abstract}

\keywords {Cohomology of Banach algebras}

\thanks{I am indebted to the Kansas State University, the University
of California at San Diego and the Mathematical Sciences Research
Institute for hospitality while this work was carried out. Research at
MSRI is supported in part by NSF grant DMS-9022140.  It was supported
also by an RFFI grant; Project 93-011-156.  }

\maketitle

\section{Introduction}

     Much interest has attached in recent years to the computation of cyclic
(co)hom\-ol\-ogy groups; see [Lo] for many references.  Most of the literature has been
devoted to the purely algebraic context, but there have also been papers
addressing the calculation of the Banach version of these groups  for Banach algebras,
and in particular $C^*$-algebras: see, for example, [ChS1], [Wo] and [He2]. There is an
effective tool for computing cyclic (co)homology: the Connes-Tsygan exact sequence,
which connects the cyclic (co)homology of many algebras with their simplicial 
(co)homology. However, it remains the case that these groups can only be calculated for
a restricted range of algebras. The purpose of this paper is to describe a technique
for the calculation of the Banach simplicial and cyclic cohomology groups of 
Banach algebras,
to establish the basic properties of this technique and to apply it to some natural
classes of algebras. The technique involves, for a Banach algebra  $A$, not 
necessarily unital, the notion of Banach cyclic cohomology 
groups $ {\cal HC}^n_B(A)$ of $A$ relative to a closed 
subalgebra $B$ of $A$. This concept was introduced and exploited  
by L. Kadison [Ka1]in the algebraic theory. We establish a connection between 
the Banach 
cyclic cohomology group $ {\cal HC}^n(A)$  of $A$ and the Banach
 $B$-relative cyclic cohomology group $ {\cal HC}^n_B(A) $ 
 of $A$. The key result (Theorem 4.1) is that, 
for a Banach algebra $A$ with a 
 bounded approximate identity and an amenable closed subalgebra $B$ of $A$,
 ${\cal HC}^n(A) = {\cal HC}^n_B(A) $ for all $n \ge 0$.
 With the aid of this theorem we show, for example, 
 that if ${\cal R}$ is a von Neumann algebra then, for each $n \ge 0$, 
${\cal HC}^n({\cal R})$ is the direct sum of the cyclic cohomologies 
$ {\cal HC}^n({\cal R}_{I_f})$ and  
$ {\cal HC}^n({\cal R}_{II_1})  $ of the finite von Neumann algebras
$ {\cal R}_{I_f}$ and  ${\cal R}_{II_1}$  of types $I$ and $II$
  appearing in the standard central direct 
summand decomposition of ${\cal R}$ (Corollary 4.5).
( In this paper equality of cohomology groups means topological isomorphism 
of seminormed spaces.)

     In establishing the connection between the Banach
cyclic cohomology groups and the Banach relative cyclic cohomology groups of a 
Banach algebra 
$A$, we need to know the connection between the Banach simplicial cohomology groups  
 ${\cal H}^n(A, A^*)$ 
of $A$ and its relative analogue $ {\cal H}^n_B(A, A^*) $. Here $A^*$ is the dual Banach
space of $A$. 
 In Section 1 we prove Theorem 1.6 that for a Banach algebra $A$, an amenable 
closed subalgebra $B$ of $A$ and a dual $A$-bimodule $M$,
 ${\cal H}^n(A,M) = {\cal H}^n_B(A,M) $ for all $n \ge 0$. 
I am grateful to B. E. Johnson for suggesting that this
result should be true. He gave a proof in the special case of a Banach algebra $A$
with an identity, an amenable closed subalgebra 
$B$ of $A$ which contains the identity of $A$ and a unital dual $A$-bimodule $M$. He
wrote that it should not be too difficult to prove the general case by known techniques 
(personal communication, 1993). As it transpired, the identical result (for the 
special case) had been proved by different methods by  L. Kadison [Ka2].
 
    E. Christensen and A. M. Sinclair in [ChS2] used the same
version of relative Hochschild cohomology group for the computation of the  
Hochschild cohomology groups of von Neumann algebras. See also [SiSm]. 

     In Section 3 we introduce the Banach relative cyclic cohomology of a Banach algebra
and show that the relative Connes-Tsygan exact sequence exists for a Banach algebra
with a bounded approximate identity. We do this using ideas of A. Ya.
Helemskii [He2]. In this section we also show the existence of morphisms of certain
 Connes-Tsygan exact sequences.

     The results of Section 1 allow us to find, in Section 2, a connection between 
the cohomology groups of a Banach algebra $A$ and that of a quotient algebra $A/I$ 
by an amenable closed bi-ideal $I$. We prove Theorem 2.1 that in this
case, for a dual $A/I$-bimodule $M$, 
 ${\cal H}^n(A,M) = {\cal H}^n(A/I,M) $ for all $n \ge 0$. Thus we can obtain,
in Theorem 2.4, the following information about the Banach simplicial cohomology groups: 
${\cal H}^n(A,A^*) = {\cal H}^n(A/I,(A/I)^*) $ for all $n \ge 2$ and amenable $I$. 
A connection between the 
Banach cyclic
cohomology groups of a Banach algebra $A$ and that of a quotient algebra $A/I$ 
by an amenable closed bi-ideal $I$ is given in Theorem 4.2 of  Section 4.

  Finally,  I would like to thank J. H. Anderson and N. J. Young for some helpful comments.

\section{Definitions and notation}

     Let $A$ be a Banach algebra, not necessarily unital, and let $A_+$ be the 
unitization of $A$. We denote by $e_+$ the adjoined identity and by $e$ an identity
of $A$ when it exists.

  We recall some notation and terminology used in the homological theory of
Banach algebras.
 Let $A$ be a Banach algebra, not necessarily unital, and 
 let $X$ be a Banach $A$-bimodule. 
We define {\it  an $n$-cochain} to be a bounded $n$-linear 
operator of $A \times \dots \times A$ into $X$ and we denote 
 the space of $n$-cochains by $ C^n(A,X)$. For $n=0$ the space  $ C^0(A,X)$ is
defined to be $X$. Let us consider the {\it standard cohomological complex}

\hspace{2.5cm}
$ 0 \rightarrow C^0(A,X) \stackrel { \delta^0} { \rightarrow} 
\dots \rightarrow C^n(A,X) \stackrel { \delta^n} { \rightarrow}   C^{n+1}(A,X) 
\rightarrow \dots,
\hfill {({\cal C}(A,X))} $
\vspace*{0.2cm}
\newline
where the coboundary operator $ \delta^n $ is defined by
\[
(\delta^n f)(a_1,...,a_{n+1}) = a_1 \cdot f(a_2,...,a_{n+1}) +
\]
\[
 \sum_{i=1}^n (-1)^i f(a_1,...,a_i a_{i+1},...,a_{n+1}) +
(-1)^{n+1} f(a_1,...,a_n) \cdot a_{n+1}.
\]
The kernel of $\delta^n$ in  $ C^n(A,X)$ is denoted by $ Z^n(A,X)$ and its elements are
called {\it $n-$cocycles}. The image of  $\delta^{n-1}$ in  $ C^n(A,X)$ $(n \ge 1)$ is
denoted by  $ N^n(A,X)$ and its elements are
called {\it $n-$coboundaries}. An easy computation yields 
$\delta^{n+1} \circ \delta^n = 0$, $n \ge 0$.

{\sc 0.1 Definition.}~~  The $n$th cohomology group of ${\cal C}(A,X)$ 
is called {\it  the $n$-dimensional Banach
cohomology group of $A$ with coefficients in $X$}. It is denoted by $ {\cal H}^n(A,X)$.

    Thus $ {\cal H}^n(A,X)  =  Z^n(A,X)/N^n(A,X); $ it is a complete seminormed space.

    Recall that a Banach $A$-bimodule $M = (M_*)^*$, where $M_*$ is a Banach
$A$-bimodule, is called {\it dual}. A Banach algebra $A$ such that 
 $ {\cal H}^1(A,M) = \{ 0 \}$
for all dual $A$-bimodules  $M$ is called {\it amenable}.

   Now let $S$ be a closed subalgebra of $A_+$. 
  We denote by $ C_S^n(A,X)$ the closed subspace of  $ C^n(A,X)$ of $n$-cochains 
$\rho$ such that
\[
 \rho (s a_1, a_2,...,a_n)  = s \cdot \rho (a_1,...,a_n) ,
\]
\[
  \rho (a_1,..., a_{i-1}, a_i s, a_{i+1}, a_{i+2},...,a_n) =
 \rho (a_1,..., a_{i-1}, a_i, s a_{i+1}, a_{i+2},...,a_n)
\]
and
\[
 \rho (a_1, a_2,...,a_n s)  = \rho (a_1,...,a_n) \cdot s
\]
for all $a_1, a_2,..., a_n \in A$, $s \in S$ and $1 \le i \le n$.
These cochains we shall call {\it $S$-relative $n$-cochains}.
 For $n=0$ the space  $ C_S^0(A,X)$ is
defined to be $ {\rm Cen}_S X  \stackrel {def} {=} \{ x \in X : 
\; s \cdot x = x \cdot s \; {\rm for \; all} \; s \in S \}$.

Note that, for each $\rho \in  C_S^n(A,X)$,  $\delta^n \rho$ is also an $S$-relative
cochain. Therefore there is a subcomplex in 
 ${\cal C}(A,X)$ formed by the spaces    $ C_S^n(A,X)$. We denote this subcomplex by
  ${\cal C}_S(A,X)$.
The kernel of $\delta^n$ in  $C_S^n(A,X)$ is denoted by $ Z_S^n(A,X)$ and 
its elements are
called {\it $S$-relative $n-$cocycles}. The image of  $\delta^{n-1}: C_S^{n-1}(A,X) 
\rightarrow C_S^n(A,X)$ $(n \ge 1)$ is denoted by  $ N_S^n(A,X)$ and its elements are
called {\it $S$-relative $n-$coboundaries}. 

{\sc 0.2 Definition.}~~  The $n$th cohomology group of ${\cal C}_S(A,X)$ 
is called  {\it the $n$-dimen\-sional Banach $S$-relative
cohomology group of $A$ with coefficients in $X$}. It is denoted by $ {\cal H}_S^n(A,X)$.

    When $S = \Bbb{C} e_+$ the subscript $S$ is unnecessary and we omit it.
 
    Throughout the paper $id$ denotes the identity operator. We denote the
projective tensor product of Banach spaces by  $\hat{\otimes}$ and the
projective tensor product of left and right Banach $A$-bimodules  by  $\hat{\otimes}_A$
[Rie].

\section{Relative cohomology of Banach algebras}

     We need a strengthening of Theorem 4.1 of [JKR] to prove the isomorphism
of the cohomology and the relative cohomology of a Banach algebra $A$ for 
dual $A$-bimodules.

{\sc 1.1 Proposition.}~~ {\it Let $A$ be a Banach algebra, let $B$ be an amenable 
closed subalgebra of $A$, let $M$ be a dual $A$-bimodule and let $n \ge 1$. Suppose
$\rho \in C^n(A,M)$  is such that 
\[
(\delta^n \rho)(a_1,...,a_{n+1}) = 0
\]
if any one of $a_1,...,a_{n+1}$ lies in $B$. Then there exists $ \xi \in 
 C^{n-1}(A,M)$ such that 
\[
(\rho - \delta^{n-1} \xi)(a_1,...,a_n) = 0
\]
if any one of $a_1,...,a_n$ lies in $B$.}

{\it The proof} is the same as that of Theorem 4.1 of [JKR].
\hfill $\Box$

     The following is essentially Lemma 4.1 of [Rin], but with a weakening of the
 hypothesis.

{\sc 1.2 Lemma.}~~ {\it Let $A$ be a Banach algebra, let $B$ be a 
closed subalgebra of $A$, let $X$ be a Banach $A$-bimodule and let $n \ge 1$.
Suppose
$\rho \in C^n(A,X)$ is such that
\[
(\delta^n \rho)(a_1,...,a_{n+1}) = 0
\]
if any one of $a_1,...,a_{n+1}$ lies in $B$ and
\[
 \rho(a_1,...,a_n) = 0
\]
if any one of $a_1,...,a_n$ lies in $B$. Then $ \rho \in C_{B}^n(A,X)$. 
 }

{\it The proof} is the same as that of Lemma 4.1 of [Rin].
\hfill $\Box$

{\sc 1.3 Corollary.}~~ {\it Let $A$ be a Banach algebra, let $B$ be an amenable 
closed subalgebra of $A$, let $M$ be a dual $A$-bimodule and let $n \ge 1$. Suppose
$\rho \in Z^n(A,M)$. 
Then there exists $ \xi \in 
 C^{n-1}(A,M)$ such that 
\[
(\rho - \delta^{n-1} \xi)(a_1,...,a_n) = 0
\]
if any one of $a_1,...,a_n$ lies in $B$. Moreover
$(\rho - \delta^{n-1} {\xi}) \in Z_{B}^n(A,M)$.}

{\sc 1.4 Proposition.}~~ {\it Let $A$ be a Banach algebra, let $B$ be a
closed subalgebra of $A$ with a bounded approximate identity $e_\nu, \nu \in
\Lambda,$ let $M = (M_*)^*$ be a dual $A$-bimodule 
and let $n \ge 1$. Suppose
$\rho \in C_B^n(A,M)$ is such that  
\begin{equation}
\label{assumption}
(\delta^n \rho)(a_1,...,a_{n+1}) = 0
\end{equation}
if any one of $a_1,...,a_{n+1}$ lies in $B$. Then there exists $ \xi_B \in 
 C^{n-1}_B(A,M)$ such that 
\[
(\rho - \delta^{n-1} \xi_B)(a_1,...,a_n) = 0
\]
if any one of $a_1,...,a_n$ lies in $B$.}

{\it Proof.} For $n=1$, by assumption, for each $b \in B$ and $\nu \in \Lambda$,
\[
0 = (\delta^1 \rho)(e_\nu,b) = 
e_\nu \cdot \rho(b) -  \rho(e_\nu b) + \rho(e_\nu) \cdot b =
 \rho(e_\nu) \cdot b,
\]
since $\rho \in C_B^1(A,M)$ and  $e_\nu \in B$.
So we obtain 
\[
 \rho(b) = \lim_{\nu}  \rho(e_\nu b) =  
 \lim_{\nu} \rho(e_\nu) \cdot b = 0.
\]
Hence we can take $\xi_B = 0$.

      For $n > 1$, we construct, inductively on $k$, $ \xi_1,...,\xi_k$ in
 $ C^{n-1}_B(A,M)$ such that 
\[
(\rho - \delta^{n-1} \xi_k)(a_1,...,a_n) = 0
\]
if any one of $a_1,...,a_k$ lies in $B$ for $1 \le k \le n$. The conclusion 
of the proposition then follows, with $\xi_B = \xi_n$.

     To construct $\xi_1$, we consider $g_{\nu} \in  C^{n-1}(A,M)$ given by
\[
 g_{\nu}(a_1,...,a_{n-1}) = \rho(e_\nu,a_1,...,a_{n-1}),
\]
for $a_1,...,a_{n-1} \in A$ and $\nu \in \Lambda$. 
Since the Banach space  $ C^{n-1}(A,M)$ is the dual space of  
$A \hat{\otimes} \dots  \hat{\otimes} A  \hat{\otimes} M_*$,
for the bounded net $g_{\nu}, \nu \in \Lambda$, there exists a subnet 
 $g_{\mu}, \mu \in \Lambda'$, which ${\rm weak^*}$ converges to some cochain 
 $g \in C^{n-1}(A,M)$. It is routine to check that  $g \in C^{n-1}_B(A,M)$.

     By assumption, for each $b \in B$ and $\nu \in \Lambda$,
\[
0 = \delta^n \rho(e_\nu,b,a_2,...,a_n) =
\]
\[
 e_{\nu} \cdot \rho(b,a_2,...,a_n) -
 \rho (e_{\nu}b,a_2,...,a_n) + \rho (e_{\nu},b a_2,...,a_n) 
\]
\begin{equation}
\label{delta}
+ \sum_{i=2}^{n-1} (-1)^{i+1}  \rho (e_{\nu},b,a_2,...,a_i a_{i+1},...,a_n) +
(-1)^{n+1} \rho (e_{\nu},b,a_2,,...,a_{n-1}) \cdot a_n = 
\end{equation}
\[
 \rho (e_{\nu},b a_2,...,a_n) 
\]
\[
+ \sum_{i=2}^{n-1} (-1)^{i+1}  \rho (e_{\nu},b,a_2,...,a_i a_{i+1},...,a_n) +
(-1)^{n+1} \rho (e_{\nu},b,a_2,,...,a_{n-1}) \cdot a_n .
\]
Thus
\[
 \rho(b,a_2,...,a_n) - \delta^{n-1} g(b,a_2,...,a_n) = 
\]
\[
\rho(b,a_2,...,a_n) -
 b \cdot g(a_2,...,a_n) + g (ba_2,...,a_n) 
\]
\[
+ \sum_{i=2}^{n-1} (-1)^{i+1} g(b,a_2,...,a_i a_{i+1},...,a_n) +
(-1)^{n+1}  g(b,a_2,,...,a_{n-1}) \cdot a_n = 
\]
\[
\rho(b,a_2,...,a_n) -
 \lim_{\mu} b \cdot \rho(e_{\mu}, a_2,...,a_n) +
 \lim_{\mu}[\rho (e_{\mu},b a_2,...,a_n) 
\]
\[
+ \sum_{i=2}^{n-1} (-1)^{i+1}  \rho (e_{\mu},b,a_2,...,a_i a_{i+1},...,a_n) +
(-1)^{n+1} \rho (e_{\mu},b,a_2,,...,a_{n-1}) \cdot a_n].
\]
The first two terms cancel, and the remaining ones add up to zero by (\ref {delta}).
This proves the existence of a suitable cochain $\xi_1 = g$.

     Suppose now that $1 \le k < n$, and a suitable cochain $\xi_k \in
 C_{B}^{n-1}(A,M)$  has been constructed. With $\rho - \delta^{n-1} 
\xi_k \in  C_{B}^n(A,M)$ denoted by $\sigma$,
\begin{equation}
\label{induct}
\sigma(a_1,...,a_n) = 0
\end{equation}
if any one of $a_1,...,a_k$ lies in $B$. In order to continue the inductive process
(and so complete the proof of the theorem), it suffices to construct $\zeta$ in 
 $C_{B}^{n-1}(A,M)$ such that $\sigma - \delta^{n-1} \zeta $
 vanishes whenever any one of its first $k+1$ arguments lies in $B$. For then we have
$ \rho - \delta^{n-1} (\xi_k + \zeta) =  \sigma - \delta^{n-1} \zeta $, and we may
take $\xi_{k+1} = \xi_k + \zeta$. To this end, we consider 
 $g_{\nu} \in  C^{n-1}(A,M)$ given by
\[
 g_{\nu}(a_1,...,a_{n-1}) = \sigma(a_1, \dots, a_k, e_\nu, a_{k+1}, \dots,a_{n-1}),
\]
for $a_1,...,a_{n-1} \in A$ and $\nu \in \Lambda$. 
By the same arguments as in the case $k=1$,  there exists a subnet 
 $g_{\mu}, \mu \in \Lambda'$, which ${\rm weak^*}$ converges to some cochain 
 $g \in C^{n-1}(A,M)$. It can be checked that  $g \in C^{n-1}_B(A,M)$ and 
\begin{equation}
\label{propertyg}
 g(a_1,...,a_{n-1}) =0
\end{equation}
if any one of $a_1,...,a_k$ lies in $B$.

     In view of the assumption (\ref{assumption}), for any $b \in B$ and 
 $\nu \in \Lambda$,
\[
\displaylines{
 \delta^n \sigma (a_1, \dots, a_k, e_\nu, b, a_{k+2}, \dots, a_n) 
\hfill\cr\hfill{}=
 (\delta^n \rho - (\delta^n \circ \delta^{n-1})(\xi_k))
(a_1, \dots, a_k, e_\nu, b, a_{k+2}, \dots, a_n) = 0.}
\]
Hence by the coboundary formula,

\[
0 = \delta^n \sigma (a_1, \dots, a_k, e_\nu, b, a_{k+2}, \dots, a_n) =
\]
\[
 a_1 \cdot \sigma(a_2, \dots, a_k, e_\nu, b, a_{k+2}, \dots, a_n)  
\]
\[
+
\sum_{i=1}^{k-1} (-1)^i 
\sigma(a_1,\dots,a_i a_{i+1},\dots,a_k,e_\nu,b,a_{k+2},\dots,a_n)  
\]
\[ 
+ (-1)^k \sigma(a_1, \dots, a_k e_\nu, b, a_{k+2}, \dots, a_n)
+ (-1)^{k+1} \sigma(a_1, \dots, a_k, e_\nu b, a_{k+2}, \dots, a_n)
\]
\[
+ (-1)^{k+2} \sigma(a_1, \dots, a_k, e_\nu, b a_{k+2}, \dots, a_n)
\]
\[
+ \sum_{i=k+2}^{n-1} (-1)^{i+1} 
\sigma(a_1, \dots, a_k, e_\nu, b, a_{k+2}, \dots, a_i a_{i+1}, \dots, a_n) 
\]
\[
+ (-1)^{n+1} \sigma(a_1, \dots, a_k, e_\nu, b, a_{k+2}, \dots, a_{n-1}) \cdot a_n .
\]
By the inductive hypothesis (\ref{induct}) the first two terms vanish, and since 
$\sigma \in C^n_B(A,M)$, the third and fourth cancel. Thus 
\[
0 = (-1)^{k+2} \sigma(a_1, \dots, a_k, e_\nu, b a_{k+2}, \dots, a_n)
\]
\begin{equation}
\label{deltak}
+ \sum_{i=k+2}^{n-1} (-1)^{i+1} 
\sigma(a_1, \dots, a_k, e_\nu, b, a_{k+2}, \dots, a_i a_{i+1}, \dots, a_n) 
\end{equation}
\[
 + (-1)^{n+1} \sigma(a_1, \dots, a_k, e_\nu, b, a_{k+2}, \dots, a_{n-1}) \cdot a_n .
\]

Now consider
\[
( \sigma - (-1)^k \delta^{n-1} g)(a_1, \dots, a_k, b, a_{k+2}, \dots, a_n) = 
\]
\[
\sigma(a_1, \dots, a_k, b, a_{k+2}, \dots, a_n) 
\]
\[
+ (-1)^{k+1} a_1 \cdot g(a_2,\dots, a_k, b, a_{k+2}, \dots, a_n) 
\]
\[
+
\sum_{i=1}^{k-1} (-1)^{i+k+1} 
g(a_1,\dots,a_i a_{i+1},\dots, a_k, b, a_{k+2}, \dots, a_n) 
\]
\[
- g (a_1, \dots, a_k b, a_{k+2}, \dots, a_n) +  
g (a_1, \dots, a_k, b a_{k+2}, \dots, a_n)
\]
\[
+ \sum_{i=k+2}^{n-1} (-1)^{i+k+1} 
 g (a_1, \dots, a_k, b, a_{k+2}, \dots, a_i a_{i+1},...,a_n) 
\]
\[
+
(-1)^{n+k+1}   g (a_1, \dots, a_k, b, a_{k+2}, \dots, a_{n-1}) \cdot a_n .
\]
By (\ref{propertyg}), the second and third terms vanish. So, by definition of $g$,
we have
\[
( \sigma - (-1)^k \delta^{n-1} g)(a_1, \dots, a_k, b, a_{k+2}, \dots, a_n) = 
\]
\[
  \sigma(a_1, \dots, a_k, b, a_{k+2}, \dots, a_n) -
 \lim_{\mu} \sigma(a_1, \dots, a_k b, e_\mu, a_{k+2}, \dots, a_n)
\]
\[
+ (-1)^k \lim_{\mu}[(-1)^{k+2} \sigma(a_1, \dots, a_k, e_\mu, b a_{k+2}, \dots, a_n)
\]
\[
+ \sum_{i=k+2}^{n-1} (-1)^{i+1} 
\sigma(a_1, \dots, a_k, e_\mu, b, a_{k+2}, \dots, a_i a_{i+1}, \dots, a_n) 
\]
\[
 + (-1)^{n+1} \sigma(a_1, \dots, a_k, e_\mu, b, a_{k+2}, \dots, a_{n-1}) \cdot a_n].
\]
The first two terms cancel, and the remaining ones add up to zero by (\ref{deltak}).
This shows that, if $\zeta = (-1)^k g$, then $\sigma -  \delta^{n-1} \zeta$ vanishes
when its $(k+1)$th argument lies in $B$. When $a_i = b \in B $ for some $i, 1 \le i \le
k$, by the inductive hypothesis (\ref{induct}) and (\ref{propertyg}), we obtain
\[
(\sigma - (-1)^k \delta^{n-1} g)(a_1, \dots, b, \dots, a_n) =
\]
\[
 (-1)^{k+i}g(a_1, \dots, a_{i-1} b, a_{i+1}, \dots, a_n) +
(-1)^{k+i+1}g(a_1, \dots, a_{i-1}, b a_{i+1}, \dots, a_n) = 0, 
\]
since $g \in C_B^{n-1}(A,M)$. Thus $\sigma -  \delta^{n-1} \zeta$ vanishes
when any of its first $k+1$ arguments lies in $B$.
 As noted above, this completes the proof of the proposition.
\hfill $\Box$

{\sc 1.5 Proposition.}~~ {\it Let $A$ be a Banach algebra, let $B$ be an amenable 
closed subalgebra of $A$, let $M$ be a dual $A$-bimodule and let $n \ge 1$. Suppose
$\rho \in C_{B}^n(A,M) \cap N^n(A,M) $. 
Then there exists $ \xi \in C_{B}^{n-1}(A,M)$ such that 
\[ 
\delta^{n-1} \xi = \rho.
\] }

{\it Proof.} For $n = 1$, by Proposition 1.4, there exists $ \xi_B \in C_{B}^0(A,M) =
{\rm Cen}_B M $ such that 
\[
(\rho - \delta^0 \xi_B)(a_1) = 0
\]
if $a_1 \in B$.
By assumption,  there exists $ \xi_1 \in C^0(A,M) = M$ such that $\rho = 
\delta^0 {\xi_1}.$ Hence 
\[
 (\rho - \delta^0 {\xi_B})(a_1) = ( \delta^0 \xi_1 - \delta^0 \xi_B)(a_1)
                    = \delta^0( {\xi_1} -  {\xi_B})(a_1)
                    = a_1.( {\xi_1} -  {\xi_B}) - ( {\xi_1} -  {\xi_B}).a_1
                    = 0
\]
if $a_1 \in B$. This implies $ {\xi_1} -  {\xi_B} \in 
{\rm Cen}_B M $ , and so  $ {\xi_1} \in {\rm Cen}_B M $.

      For $n \ge 2$, by Proposition 1.4, there exists $ \xi_B \in C_{B}^{n-1}(A,M)$ 
such that 
\[
(\rho - \delta^{n-1} \xi_B)(a_1,...,a_n) = 0
\]
if any one of $a_1,...,a_n$ lies in $B$,
and by Lemma 1.2, $\rho_B \stackrel {def} {=} (\rho - \delta^{n-1} {\xi_B}) 
\in C_{B}^n(A,M)$.
By assumption,  there exists $ \xi_1 \in C^{n-1}(A,M)$ such that 
$\rho = \delta^{n-1} {\xi_1}.$ Hence 
\[
 \rho_B =\rho - \delta^{n-1} {\xi_B} = \delta^{n-1} \xi_1 - \delta^{n-1} \xi_B
                    = \delta^{n-1}( {\xi_1} -  {\xi_B}).
\]
Further, $\eta  \stackrel {def} {=} {\xi_1} -  {\xi_B}$ satisfies the assumption of
Proposition 1.1, and so there exists $ \beta 
 \in C^{n-2}(A,M)$ such that 
\[
 \eta_B(a_1,...,a_{n-1}) \stackrel {def} {=} 
(\eta - \delta^{n-2} \beta)(a_1,...,a_{n-1}) = 0
\]
if any one of $a_1,...,a_{n-1}$ lies in $B$.
     Therefore 
\[
 \rho_B  = \delta^{n-1}( {\xi_1} -  {\xi_B})
         = \delta^{n-1} \eta
         =  \delta^{n-1}(\eta_B +  \delta^{n-2} \beta)
         =  \delta^{n-1} \eta_B.
\]
By Lemma 1.2, $\eta_B \in  C_{B}^{n-1}(A,M)$; this implies
\[
 \rho = \rho_B + \delta^{n-1} {\xi_B}
       = \delta^{n-1} \eta_B + \delta^{n-1} {\xi_B}
       = \delta^{n-1}(\eta_B + \xi_B)
\]
and so, for 
  $ \xi = \eta_B + \xi_B \in C_{B}^{n-1}(A,M)$ ,
we have $\delta^{n-1} \xi = \rho.$
\hfill $\Box$

     As was said in the introduction, the following result is based on a communication
of B. E. Johnson.

{\sc 1.6 Theorem.}~~ {\it Let $A$ be a Banach algebra, let $B$ be an amenable 
closed subalgebra of $A$ and let $M$ be a dual $A$-bimodule.
Then
\[
{\cal H}^n(A,M) = {\cal H}_B^n(A,M)
\] 
for all $n \ge 0$.}

{\it Proof.} The inclusion morphism of cochain objects
$ C_B^n(A,M)  \rightarrow  C^n(A,M) $
induces  a morphism of complexes $ {\cal C}_B(A,M)  \rightarrow   {\cal C}(A,M)$ and
hence morphisms 
\[
    {\cal F}_n: {\cal H}_B^n(A,M)  \rightarrow   {\cal H}^n(A,M)
\]
given, for each $\rho \in Z_{B}^n(A,M)$,  by 
  $  {\cal F}_n(\rho + N_{B}^n(A,M)) = \rho + N^n(A,M)$ (see, for example,
[He1; Section 0.5.3]).

    For $n = 0$ we have $ {\cal H}_B^0(A,M) = {\cal H}^0(A,M) =
{\rm Cen}_A(M) $.
      
      In the case where $n \ge 1$, the morphism $ {\cal F}_n$ is injective by
Proposition 1.5 and surjective by Corollary 1.3. Hence, by Lemma 0.5.9 [He1], 
 $ {\cal F}_n$ is a topological isomorphism. 
\hfill $\Box$

     In the particular case when $A$ is a unital $C^*$-algebra, 
$B$ is the $C^*$-algebra generated by an amenable group of unitaries
and $M$ is a dual normal $A$-bimodule, the statement of Theorem 1.6 is given in [SiSm];
see Theorem 3.2.7 [SiSm].

{\sc 1.7 Proposition.}~~ {\it Let $A_i$ be a Banach algebra with 
identity $e_i, i = 1, ..., m,$ let $A$ be the Banach algebra direct sum 
$ \bigoplus_{i=1}^m A_i$ and let $M$ be a dual $A$-bimodule.
Then the canonical projections from $A$ to $A_i, i = 1, ..., m,$ induce a 
topological isomorphism
of Banach spaces $C_B^n(A,M)  \rightarrow  \bigoplus_{i=1}^m C^n(A_i, e_i M e_i)$, where
 $ n \ge 0$ and $B$ is the Banach subalgebra of $A$ generated by 
$ \{e_i, i = 1, ..., m\}$. Hence
\[ 
{\cal H}^n(A,M) =  \bigoplus_{i=1}^m {\cal H}^n(A_i, e_i M e_i) 
\]
for all $n \ge 0$.}

Here $e_i M e_i = \{e_i \cdot x \cdot e_i, x \in M\}$ is a Banach $A$-bimodule.

{\it Proof.} For $n = 0$ it can be checked that 
\[ 
C_B^0(A,M) = {\rm Cen}_B(M) = \bigoplus_{i=1}^m (e_i M e_i) =
 \bigoplus_{i=1}^m C^0(A_i, e_i M e_i).
\]
 For $n \ge 1$ and $\rho \in C_B^n(A,M)$, we have
\[
\rho(a_1,...,a_n) = \rho( \sum_{i=1}^m a_1 e_i, a_2,...,a_n) = 
   \sum_{i=1}^m e_i \cdot \rho(a_1 e_i, a_2 e_i,...,a_n e_i) \cdot e_i
\]
for $a_1, \dots, a_n \in A$.
We define cochain maps from $ C_B^n(A,M)$ to       
$ \bigoplus_{i=1}^m C^n(A_i,e_i M e_i)$ and back by
\[
    {\cal J}_n: \rho \mapsto   (\rho_1, \dots, \rho_m)
\]
where $\rho_i(a^i_1, \dots, a^i_n) = 
\rho(a^i_1, \dots, a^i_n)$ 
for $a^i_1, \dots, a^i_n \in A_i$ and for $i = 1, \dots, m,$
and 
\[
    {\cal G}_n:   (\rho_1, \dots, \rho_m) \mapsto   \rho
\]
is given by $\rho(a_1,...,a_n) =  
   \sum_{i=1}^m \rho_i(a_1 e_i, a_2 e_i,...,a_n e_i) $
for $a_1, \dots, a_n \in A$.

     It is clear that $ {\cal J}_n \circ  {\cal G}_n = id$, 
 $ {\cal G}_n \circ  {\cal J}_n = id$ and maps  $ {\cal J}_n$ and $ {\cal G}_n $ are
bounded. Thus there is a topological isomorhism of complexes
${\cal C}_B(A,M)  \rightarrow  \bigoplus_{i=1}^m {\cal C}(A_i, e_i M e_i)$, and so
\[ 
{\cal H}_B^n(A,M) =  \bigoplus_{i=1}^m {\cal H}^n(A_i, e_i M e_i). 
\]
Note that $ B = \bigoplus_{i=1}^m \Bbb{C} e_i$ is amenable. Hence, 
by Theorem 1.6, ${\cal H}^n(A,M) ={\cal H}_B^n(A,M) $ for all $n \ge 0$.
The result now follows directly. 
\hfill $\Box$

{\sc 1.8 Proposition.}~~ {\it Let ${\cal R}$ be a von Neumann algebra, let
\[
{\cal R} = {\cal R}_{I_f}  \oplus {\cal R}_{I_{\infty}}  \oplus {\cal R}_{II_1} 
\oplus {\cal R}_{II_\infty} \oplus {\cal R}_{III}
\]
be the central direct summand decomposition of ${\cal R}$ into von Neumann algebras 
of types 
$I_f (finite) , I_{\infty}, II_1, II_{\infty}, III$  with the identity $e$ of 
${\cal R}$ decomposing as 
$e =  e_{I_f} \oplus e_{I_\infty} \oplus e_{II_1} \oplus e_{II_\infty} \oplus e_{III}$ 
([Sa; Section 2.2]), and let $M$ be a dual ${\cal R}$-bimodule.
Then 
 
{\rm (i)}
\[ 
{\cal H}^n({\cal R},M) =
 {\cal H}^n({\cal R}_{I_f}, e_{I_f} M e_{I_f})  \oplus
{\cal H}^n({\cal R}_{I_\infty}, e_{I_\infty} M e_{I_\infty})  \oplus
\]
\[
 {\cal H}^n({\cal R}_{II_1}, e_{II_1} M e_{II_1})  \oplus 
{\cal H}^n({\cal R}_{II_\infty}, e_{II_\infty} M e_{II_\infty})  \oplus
 {\cal H}^n({\cal R}_{III}, e_{III} M e_{III})  
\]
for all $n \ge 0$;

{\rm (ii)} in particular, 
\[
{\cal H}^n({\cal R},{\cal R}^*) =  {\cal H}^n({\cal R}_{I_f},{\cal R}_{I_f}^*) 
\oplus {\cal H}^n({\cal R}_{II_1},{\cal R}_{II_1}^*) 
\]
for all $n \ge 0$.}

{\it Proof.} Part (i) follows from Proposition 1.7. In part (ii) we apply (i) to
$M = {\cal R}^* $. By Proposition 2.2.4 [Sa], for a properly infinite von Neumann
algebra ${\cal U}$ there exists a sequence $(p_m)$ of mutually orthogonal, equivalent
projections in  ${\cal U}$ with $p_m \sim e$. Thus, by Theorem 2.1 [Fa], each hermitian
element of  ${\cal U}$ is the sum of five commutators. Hence there are no 
non-zero bounded traces on ${\cal U}$. Thus, by virtue of Corollary 3.3 
[ChS1], for a von Neumann algebra ${\cal U}$ of one of the types $I_{\infty}, 
II_{\infty}$ 
or $III$, the simplicial cohomology groups  
${\cal H}^n({\cal U},{\cal U}^*) =\{ 0 \}$ for all $n \ge 0$.
 \hfill $\Box$

{\sc 1.9 Remark.}~~ As for finite von Neumann algebras of type $I$, 
by [Sa; Theorem 2.3.2], they are the $l_{\infty}$-direct sum of type $I_m$ von Neumann 
algebras ${\cal R}_{I_m}$, where $m < \infty$. By Theorem 2.3.3 [Sa] and 
results of
Section 1.22 [Sa], ${\cal R}_{I_m}$ is $*$-isomorphic to the $C^*$-tensor product 
$Z \otimes_{min} {\cal B}(H)$, where  $Z$ is the centre of ${\cal R}_{I_m}$ and 
$\dim (H) = m$.
Hence, by Theorem 7.9 [Jo1], finite von Neumann algebras ${\cal R}_{I_m}$ of type $I_m$ 
are amenable, 
and so their simplicial cohomology groups vanish 
${\cal H}^n({\cal R}_{I_m},{\cal R}_{I_m}^*) = \{ 0 \}$ for all $n \ge 1$. It is still not 
clear to the author whether
${\cal H}^n({\cal R}_{I_f},{\cal R}_{I_f}^*) = \{ 0 \}$ for all $n \ge 2$.

   Note that the statement of Proposition 1.8 (i) is proved in [SiSm], 
Corollary 3.3.8, for the case of a dual normal ${\cal R}$-bimodule $M$.
However, this result does not apply to $M = {\cal R}^* $, since ${\cal R}^* $ is not in
general a normal dual module.

    Now let us consider two unital Banach algebras $A_1$ and $A_2$, 
a unital Banach $A_1-A_2$-bimodule $Y$, and the natural triangular matrix algebra  
\[
{\cal U} = \left[ \begin{array}{cc}
A_1 & Y
\\
0 & A_2
\end{array} \right]
\]
 with matrix multiplication and norm
\[
\| \left[ \begin{array}{cc}
r_1 & y
\\
0 & r_2
\end{array} \right] \| = \sup_{\| d_1\| \le 1, \| d_2\| \le 1} 
\max \{\|r_1 d_1 + y \cdot d_2\|, \|r_2 d_2\| \},
\]
 where $d_1 \in A_1, d_2 \in A_2.$ 
Let $e_{ii}, i = 1, 2,$ denote the idempotents 
$e_{11} = \left[ \begin{array}{cc}
e_{A_1} & 0
\\
0 & 0
\end{array} \right]$
and $e_{22} = \left[ \begin{array}{cc}
0 & 0
\\
0 & e_{A_2}
\end{array} \right]$ respectively. Let us also consider
the Banach subalgebra $B$ of ${\cal U}$ generated by 
$ \{e_{ii}, i = 1, 2 \}$.

{\sc 1.10 Proposition.}~~ {\it Let $A_1$ and $A_2$ be unital Banach algebras, 
 let $Y$ be a unital Banach $A_1-A_2$-bimodule, and let 
\[
{\cal U} = \left[ \begin{array}{cc}
A_1 & Y
\\
0 & A_2
\end{array} \right]
\]
be the natural triangular matrix algebra.
 Then the two canonical projections from ${\cal U}$ to $A_1$ and $A_2$ induce a
topological isomorphism of Banach spaces $C_B^n({\cal U},{\cal U}^*)  \rightarrow 
 \bigoplus_{i=1}^2 C^n(A_i, A_i^*), n \ge 0$, and hence 
\[
{\cal H}^n({\cal U},{\cal U}^*) = 
{\cal H}^n(A_1,A_1^*)  \oplus {\cal H}^n(A_2,A_2^*)
\] 
for all $n \ge 0$.} 

    The same assertion in a purely algebraic context was proved by 
L. Kadison in [Ka1]. 
 A key step in Kadison's proof is the equality  
 ${\cal H}^n({\cal U},{\cal U}^*) =
{\cal H}_B^n({\cal U},{\cal U}^*) $ for all $n \ge 0$, where $B$ is the Banach algebra
generated by $e_{11}$ and $e_{22}$. This equality remains valid in the present 
context by Theorem 1.6. The rest of the proof is as in [Ka1], with insignificant changes.

\section{The connection between the cohomologies of $A$ and $A/I$}

     Recall Proposition 5.1 [Jo1] that a quotient algebra of an amenable algebra is 
amenable, and that an extension of an amenable algebra by an
amenable  bi-ideal is an amenable algebra (this can also be found in Corollary 35 and  
Proposition 39 [He3]). The following theorem gives some additional
information about the cohomology of Banach algebras $A$ and $A/I$ without the assumption
that $A$ be amenable.

{\sc 2.1 Theorem.}~~ {\it Let $A$ be a Banach algebra and let $I$ be a 
closed two-sided ideal of $A$. Suppose that $I$ is an amenable  Banach algebra and 
$M$ is a dual $A/I$-bimodule.
Then
\[
{\cal H}^n(A,M) = {\cal H}^n(A/I,M)
\] 
for all $n \ge 0$.}

{\it Proof.} For $n = 0$ we have $ {\cal H}^0(A,M) = {\cal H}^0(A/I,M) =
{\rm Cen}_{A/I}(M) $. In the case where $n \ge 1$  the inclusion morphism of 
cochain objects
$  C^n(A/I,M)  \rightarrow C^n(A,M)$
induces a morphism of complexes $ {\cal C}(A/I,M)  \rightarrow {\cal C}(A,M)$ 
and hence morphisms 
\[
    {\cal L}_n: {\cal H}^n(A/I,M)  \rightarrow   {\cal H}^n(A,M)
\]
given, for each $\hat\rho \in Z^n(A/I,M)$,  by
\[
    {\cal L}_n ( \hat\rho + N^n(A/I,M)) =  {\rho} + N^n(A,M),
\]
where $ {\rho}(a_1,...,a_n) = \hat\rho(\theta(a_1),...,\theta(a_n))$ and
$ \theta: A \rightarrow A/I $ is the natural epimorphism. It is straight-forward to 
check that if 
 $\hat\rho = \delta^{n-1} \hat\xi$ for some
$\hat\xi \in C^{n-1}(A/I,M)$ then $\rho = \delta^{n-1} \xi$ where
 $ {\xi}(a_1,...,a_{n-1}) = \hat\xi(\theta(a_1),...,\theta(a_{n-1}))$.

    By Corollary 1.3, for $\eta \in Z^n(A,M)$,  there exists $ \xi \in 
 C^{n-1}(A,M)$ such that 
\[
(\eta - \delta^{n-1} \xi)(a_1,...,a_n) = 0
\]
if any one of $a_1,...,a_n$ lies in $I$. We can therefore define
\[
 \eta_I   \stackrel {def} {=}\eta - \delta^{n-1} \xi ~{\rm and}~
  \hat\eta_I(a_1+I,...,a_n+I) \stackrel {def} {=} \eta_I (a_1,...,a_n) .
\]
Hence for each $\eta \in Z^n(A,M)$  there exists $\hat\eta_I \in Z^n(A/I,M)$
such that
\[
    {\cal L}_n ( \hat\eta_I + N^n(A/I,M)) =  {\eta} + N^n(A,M), 
\]
and so $ {\cal L}_n $ is surjective.

     Let $ {\cal L}_n ( \hat\rho + N^n(A/I,M)) = 0$, that is,
 $ {\rho}(a_1,...,a_n) = \hat\rho(\theta(a_1),...,\theta(a_n))  \in N^n(A,M)$.
This implies that there is $\beta  \in C^{n-1}(A,M)$ such that
 $\rho = \delta^{n-1} \beta$. Further, $\beta$ satisfies the assumption of 
Proposition 1.1, and so there exists $ \alpha  \in C^{n-2}(A,M)$ such that
\[
(\beta - \delta^{n-2} \alpha)(a_1,...,a_{n-1}) = 0
\]
if any one of $a_1,...,a_{n-1}$ lies in $I$.
We can define 
\[
 \beta_I  \stackrel {def} {=} \beta - \delta^{n-2} \alpha
\]
and see that  $\rho = \delta^{n-1} \beta = \delta^{n-1} \beta_I$.
Therefore
 $\hat\rho = \delta^{n-1} \hat\beta_I$, where
\[
  \hat\beta_I(a_1+I,...,a_n+I) =\beta_I (a_1,...,a_n) .
\]
This proves the injectivity of  $ {\cal L}_n $. Hence, by Lemma 0.5.9 [He1], 
 $ {\cal L}_n$ is a topological isomorphism. 
\hfill $\Box$

{\sc 2.2 Proposition.}~~ {\it Let $A$ be a Banach algebra and let $I$ be a 
closed two-sided ideal of $A$. Suppose that $I$ has a bounded approximate identity. 
 Then ${\cal H}_I^0(A,I^*) =   {\rm Cen}_I I^* $ and
\[
{\cal H}_I^n(A,I^*) = \{ 0 \}
\] 
for all $n \ge 1$.}

{\it Proof.} Let us consider the Banach space $C_I^n(A,I^*)$
which is isometrically isomorphic to the Banach space  
$ _I h_I (A \hat{\otimes}_I  \dots  \hat{\otimes}_I A, I^*)$ of all Banach $I$-bimodule 
 morphisms from $A \hat{\otimes}_I  \dots  \hat{\otimes}_I A$  into $I^*$. 
The latter Banach space is isometricaly
isomorphic to 
$ {\rm Cen}_I (A \hat{\otimes}_I \dots \hat{\otimes}_I A \hat{\otimes}_I I)^*$
by Proposition VII.2.17 [He1].
By virtue of the assumption, $I$ has a bounded approximate identity, and
so Proposition II.3.13 [He1] gives us an isomophism of Banach $I$-bimodules
$A \hat{\otimes}_I \dots \hat{\otimes}_I A \hat{\otimes}_I I = I$.
 Therefore, there exists an isometric isomorphism of Banach spaces
\[ 
 {\cal F}_n: C_I^n(A,I^*) \rightarrow   {\rm Cen}_I I^*
\] 
for all $n \ge 0$.

       Now it is routine to check that the diagram
\[
\arraycolsep=1pt
\begin{array}{ccccccccccccc}
0 & \rightarrow & C_I^0(A,I^*) & \stackrel {\delta^0} {\rightarrow} & 
C_I^1(A,I^*) & \stackrel {\delta^1} {\rightarrow} & \dots 
& \rightarrow & C_I^n(A,I^*) &  \stackrel {\delta^n} {\rightarrow} &
 C_I^{n+1}(A,I^*) & \rightarrow & \dots 
\\
~ & ~ & \downarrow {\cal F}_0 & ~ &  \downarrow {\cal F}_1 & ~ & ~ & ~ &
 \downarrow {\cal F}_{n} & ~ &  \downarrow {\cal F}_{n +1} & ~ & ~
\\
0 & \rightarrow & {\rm Cen}_I I^* & \stackrel {\eta^0} {\rightarrow} & {\rm Cen}_I I^* &
 \stackrel {\eta^1} {\rightarrow} & \dots &
 \rightarrow & {\rm Cen}_I I^* & \stackrel {\eta^n} {\rightarrow}
 & {\rm Cen}_I I^* & \rightarrow & \dots ,
\end{array}
\] 
is commutative, where $\eta^n(f) = 0$ for all even $n$ and $\eta^n(f) = f$ for all 
odd $n$. The cohomology of the upper complex is, by definition,
 $ {\cal H}_I^n(A,I^*)$. Thus the result now follows directly.
\hfill $\Box$

{\sc 2.3 Corollary.}~~ {\it Let $A$ be a Banach algebra and let $I$ be a 
closed two-sided ideal of $A$. Suppose that $I$ is an amenable Banach algebra. 
 Then ${\cal H}^0(A,I^*) =   {\rm Cen}_I I^* $ and
\[
{\cal H}^n(A,I^*) = \{ 0 \}
\] 
for all $n \ge 1$.} 

{\it Proof.} By Theorem 1.6,
\[
{\cal H}^n(A,I^*) = {\cal H}_I^n(A,I^*) 
\] 
for all $n \ge 0$. Hence the result follows from Proposition 2.2.
 \hfill $\Box$

{\sc 2.4 Theorem.}~~ {\it Let $A$ be a Banach algebra and let $I$ be a 
closed two-sided ideal of $A$. Suppose that $I$ is an amenable Banach algebra. 
Then

   ${\rm (i)}$ 
\[
{\cal H}^n(A,A^*) = {\cal H}^n(A/I,(A/I)^*) 
\] 
for all $n \ge 2$, and the natural map from $ {\cal H}^1(A/I,(A/I)^*)$ into 
${\cal H}^1(A,A^*)$ is surjective;

    ${\rm (ii)}$ if ${\rm Cen}_I I^* = \{ 0 \} $ 
then
\[
{\cal H}^n(A,A^*) = {\cal H}^n(A/I,(A/I)^*) 
\] 
for all $n \ge 0$.} 

{\it Proof.} We consider the short exact sequence of Banach $A$-bimodules

\hspace{3.5cm}
$  0 \leftarrow A/I  \stackrel {j} {\leftarrow}    A  \stackrel {i} {\leftarrow}  I 
 \leftarrow  0, \hfill {({\cal I})} $
\newline
where $ i$ and $j$ are the natural embedding and quotient mapping 
 respectively, and its dual complex
 
\hspace{3.5cm}
$ 0 \rightarrow (A/I)^*  \stackrel {j^*} { \rightarrow} A^* \stackrel {i^*} {\rightarrow}
I^*   \rightarrow 0. \hfill {({\cal I}^*)}$
 
By virtue of its amenability, $ I$ has a bounded approximate identity and so the complex 
$~({\cal I}^*)$ is admissible. Hence, by Corollary III.4.11 [He1], there exists a long
exact sequence
\[
0 \rightarrow {\cal H}^0(A,(A/I)^*)  \rightarrow {\cal H}^0(A,A^*) \rightarrow
 {\cal H}^0(A,I^*) \rightarrow
 {\cal H}^1(A,(A/I)^*) \rightarrow {\cal H}^1(A,A^*)  \rightarrow 
\]
\vspace{0.01cm}
\[
 {\cal H}^1(A,I^*)
 \rightarrow {\cal H}^2(A,(A/I)^*)  \rightarrow {\cal H}^2(A,A^*)
 \rightarrow {\cal H}^2(A,I^*) \rightarrow  \dots
\]
\vspace{0.01cm}
\[
 \dots \rightarrow
{\cal H}^{n-1}(A,I^*) \rightarrow  {\cal H}^n(A,(A/I)^*) \rightarrow
 {\cal H}^n(A,A^*) \rightarrow {\cal H}^n(A,I^*) \rightarrow \dots .
\]
     Recall that, by Corollary 2.3,
${\cal H}^n(A,I^*) = \{ 0 \}$ for all $n \ge 1$. Thus
${\cal H}^n(A,(A/I)^*) = {\cal H}^n(A,A^*)$ (see Lemma 0.5.9 [He1])
  for all $n \ge 2$. Therefore, by Theorem 2.1,
\[
{\cal H}^n(A/I,(A/I)^*) = {\cal H}^n(A,(A/I)^*) = {\cal H}^n(A,A^*) 
\]
 for all $n \ge 2$.
\hfill $\Box$

      Note that  ${\rm Cen}_A A^* $ coincides with the space 
$$A^{tr} = \{ f \in A^* : f(ab) = f(ba) \; {\rm for \; all} \; a, b \in A \}$$ 
of continuous traces on $A$. 

     Theorem 2.4 applies whenever $I$ is a nuclear $C^*$-algebra. 
 Other examples are given by the Banach algebra $A = {\cal B}(E)$ of 
all bounded operators on a
Banach space $E$ with the property $(\Bbb{A})$, which was defined in [GJW], 
and the closed ideal $I = {\cal K}(E)$ of compact operators on $E$. In this case 
$ {\cal K}(E)$ is amenable  [GJW]. The property 
$(\Bbb{A})$ implies that ${\cal K}(E)$ contains a bounded sequence of projections of
unbounded finite rank, and from this it is easy to show (via embedding of matrix
algebras) that there is no non-zero bounded trace on  ${\cal K}(E)$.
Thus we can see from Theorem 2.4 (ii) that, for  a
Banach space $E$ with the property $(\Bbb{A})$, 
${\cal H}^n({\cal B}(E),{\cal B}(E)^*) = 
{\cal HC}^n({\cal B}(E)/{\cal K}(E),({\cal B}(E)/{\cal K}(E))^*)$  for all $n \ge 0$.
 Several classes of Banach spaces have the property $(\Bbb{A})$: 
$l_p; 1 < p < \infty; C(K)$,
where $K$ is a compact Hausdorff space; $L_p (\Omega, \mu); 1 < p < \infty$,
where $ (\Omega, \mu)$ is a measure space (for details and more examples see [Jo1] and 
[GJW]).

  In the case of $C^*$-algebras, we know that the Banach simplicial cohomology groups
vanish for $C^*$-algebras without non-zero bounded traces [ChS1; Corollary 3.3].
Therefore, for an infinite-dimensional Hilbert space $H$, we obtain
 $${\cal H}^n({\cal B}(H)/{\cal K}(H),({\cal B}(H)/{\cal K}(H))^*) = 
{\cal H}^n({\cal B}(H),{\cal B}(H)^*) = \{ 0 \}$$
for all $n \ge 0$, since   ${\cal K}(H)^{tr} = \{ 0 \}$ by [An; Theorem 2] and 
${\cal B}(H)^{tr} = \{ 0 \}$ by [Hal]. One can also see directly that the
Calkin algebra has no non-zero bounded trace, and hence has trivial Banach simplicial 
cohomology.

     Notice that, for every Banach algebra $A$, the vanishing of 
${\cal H}^n(A,A^*)$  for all $n \ge 0$ is equivalent to the vanishing of 
the Banach simplicial homology groups  ${\cal H}_n(A,A)$ for  all $n \ge 0$ 
[He1, Proposition 2.5.28].
The latter relation was established by M. Wodzicki [Wo] in the case 
  $A = {\cal B}(H)$. He also proved the 
vanishing of the (Banach) simplicial homology groups for stable $C^*$-algebras, 
i.e., for algebras isomorphic to their $C^*$-tensor product with ${\cal K}(H)$ 
for a separable $H$, and the vanishing of the algebraic simplicial homology groups for
$A = {\cal B}(H)/{\cal K}(H)$ and a separable $H$. 
 
    Recall from [Di; Sections 4.2, 4.3] that a $C^*$-algebra $A$ is called CCR (or
{\it liminary}) if $\pi(A) = {\cal K}(H)$ for each irreducible representation $(\pi, H)$
of $A$. A $C^*$-algebra $A$ is called GCR (or {\it postliminary}) if each non-zero
quotient of $A$ has a non-zero closed two-sided CCR-ideal. Finally we say that $A$ is
NGCR (or {\it antiliminary}) if it contains no non-zero closed two-sided CCR-ideal. By 
[Di; Proposition 4.3.3 and 4.3.6], each $C^*$-algebra $A$ has a largest closed
two-sided GCR-ideal $I_{\alpha}$, and $A/I_{\alpha}$ is NGCR. The following result allow
us to reduce the computation of the simplicial cohomology groups of $C^*$-algebras  
to the case of NGCR-algebras.

{\sc 2.5 Proposition.}~~ {\it Let $A$ be a $C^*$-algebra. Then 
\[
{\cal H}^n(A,A^*) = {\cal H}^n(A/I_{\alpha},(A/I_{\alpha})^*) 
\] 
for all $n \ge 1$.}

{\it Proof.} By [Ha; Corollary 4.2], ${\cal H}^1(A,A^*) = \{ 0 \}$ for every $C^*$-algebra
$A$. By Theorem 7.9 [Jo1],  $I_{\alpha}$ is amenable. Thus the result
directly follows from Theorem 2.4 (i).
\hfill $\Box$

\section{The existence of the Connes-Tsygan exact sequence}

    Let $A$ be a Banach algebra, not necessarily unital, and let $S$ be a closed 
subalgebra of $A_+$. In this section we introduce the Banach version of the concept
of $S$-relative cyclic cohomology ${\cal HC}_S^n(A)$ (compare with [Ka1]). We also show 
that the $S$-relative Connes-Tsygan exact sequence
exists for every Banach algebra $A$ with a bounded approximate identity.
This is accomplished with the aid of ideas from [He2].

    When $S = \Bbb{C} e_+$ the subscript $S$ is 
unnecessary and we omit it.

     We denote by  $ C_S^n(A), n = 0, 1,...,$ the Banach space of continuous
$(n+1)$-linear functionals on $A$ such that
\[
f(sa_0, a_1,..., a_n) = f(a_0, a_1,..., a_ns)
\]
and, for $j = 0, 1,..., n-1,$
\[
f(a_0,..., a_js, a_{j+1},..., a_n) = f(a_0,..., a_j, sa_{j+1},..., a_n)
\]
for all  $ s \in S$ and $a_0, ..., a_n \in A$; these functionals we shall call 
{\it $n$-dimensional $S$-relative cochains}. We let
\[
 t_n: C_S^n(A) \rightarrow  C_S^n(A), \qquad n = 0, 1,...,
\]
denote the operator given by
\[
t_nf(a_0, a_1,..., a_n) = (-1)^n f(a_1,..., a_n, a_0),
\]
and we set $t_0 = id$. The important point is that $t_nf$ is an $S$-relative 
cochain since
\[
\thickmuskip=.5 \thickmuskip
\medmuskip=.5 \medmuskip
t_nf(sa_0, a_1,..., a_n) = (-1)^n f(a_1,..., a_n, sa_0) =
(-1)^nf(a_1,..., a_ns, a_0) = t_nf(a_0,..., a_ns)
\]
and the other identities follow just as readily.
A cochain $f \in  C_S^n(A)$ satisfying 
$t_nf = f$ is called {\it cyclic}. We let $ CC_S^n(A)$ denote the closed subspace of
$ C_S^n(A)$ formed by the cyclic cochains.
In particular, 
\[
 CC_S^0(A) = C_S^0(A) = {\rm Cen}_S A^* = \{f \in A^*: f(sa) = f(as) \; {\rm for \; all}
\; a \in A, s \in S \}.
\]
     From the $S$-relative cochains we form the standard cohomology complex 
$\tilde{{\cal C}}_S(A)$:
\[
0 \rightarrow C_S^0(A) \stackrel {\delta^0} {\rightarrow} \dots 
\rightarrow C_S^n(A) 
 \stackrel {\delta^n} {\rightarrow}  C_S^{n+1}(A) \rightarrow \dots , 
\]
where the continuous operator $\delta^n$ is given by the formula
\[
(\delta^n f)(a_0, a_1,..., a_{n+1}) = 
\]
\[
 \sum_{i=0}^n (-1)^i f(a_0,...,a_i a_{i+1},...,a_{n+1}) +
(-1)^{n+1} f(a_{n+1} a_0,...,a_n).
\]
One can easily check that $ \delta^{n+1} \circ \delta^n$ is indeed $0$ for all $n$
and that each $\delta^nf$ is again an $S$-relative cochain. It is not difficult
to verify that every  $\delta^n$ sends a cyclic cochain again to a cyclic one.
Therefore there is a subcomplex in $\tilde{{\cal C}}_S(A)$ formed by the spaces 
 $ CC_S^n(A)$. We denote this subcomplex by  $\tilde{{\cal CC}}_S(A)$, and its 
differentials are denoted by
\[
\delta c^n: CC_S^n(A) \rightarrow  CC_S^{n+1}(A).
\]

   Note that the complex  $\tilde{{\cal C}}_S(A)$ is a subcomplex of 
 $\tilde{{\cal C}}(A)$ and  $\tilde{{\cal CC}}_S(A)$ is a subcomplex of 
 $\tilde{{\cal CC}}(A)$ respectively.

{\sc 3.1 Definition.}~~ The $n$th cohomology of $\tilde{{\cal C}}_S(A)$, denoted
by  ${\cal H}^n_S(A)$, is called the {\it $n$-dimensional Banach $S$-relative 
simplicial,} or
{\it Hochschild, cohomology group} of the Banach algebra $A$. The $n$th
 cohomology of $\tilde{{\cal CC}}_S(A)$, denoted
by   ${\cal HC}^n_S(A)$, is called the {\it $n$-dimensional Banach $S$-relative cyclic 
cohomology group} of $A$.

   Note that, by definition, $\delta c^0 = \delta^0$, so that  ${\cal HC}^0_S(A) = 
{\cal H}^0_S(A)$ coincides with the space $A^{tr} =
\{ f \in A^* : f(ab) = f(ba) \; {\rm for \; all} \; a, b \in A \}$.

   We define  ${\cal HC}^{-1}_S(A)$ to be $\{ 0 \}$.

{\sc 3.2 Remark.}~~ The canonical identification of $(n +1)-$linear functionals on $A$
 and $n-$linear operators from $A$ to $A^*$ shows that   ${\cal H}^n_S(A)$ is just 
another way of writing   ${\cal H}^n_S(A, A^*)$.

     Further, we need the following complex  
$\tilde{{\cal C}}{\cal R}_S(A)$:
\[
0 \rightarrow C_S^0(A) \stackrel {\delta r^0} {\rightarrow} \dots 
\rightarrow C_S^n(A) 
 \stackrel {\delta r^n} {\rightarrow}  C_S^{n+1}(A) \rightarrow \dots , 
\]
where the continuous operator $\delta r^n$ is given by the formula
\[
(\delta r^n f)(a_0, a_1,..., a_{n+1}) = 
 \sum_{i=0}^n (-1)^i f(a_0,...,a_i a_{i+1},...,a_{n+1}).
\]

  The $n$th cohomology of $\tilde{{\cal C}}{\cal R}_S(A)$ is denoted
by  ${\cal HR}^n_S(A)$.

     Following [He2] we consider the sequence
\[
0 \rightarrow \tilde{{\cal CC}}_S(A) \stackrel {\bar{i}} {\rightarrow} 
\tilde{{\cal C}}_S(A) \stackrel {\bar{M}} {\rightarrow} 
\tilde{{\cal C}}{\cal R}_S(A) \stackrel {\bar{N}} {\rightarrow} 
\tilde{{\cal CC}}_S(A) \rightarrow 0 
\]
of complexes in the category of Banach spaces and continuous operators, where 
 ${\bar{i}}$ denotes the natural inclusion 
\[
i_n: CC_S^n(A) \rightarrow  C_S^n(A),
\]
\[
M_n = id - t_n: C_S^n(A) \rightarrow  C_S^n(A)
\]
and
\[
N_n = id + t_n + \dots + t_n^n: C_S^n(A) \rightarrow  CC_S^n(A).
\]

{\sc 3.3 Proposition.}~~ {\it Let $A$ be a Banach algebra and let $S$ be a closed 
subalgebra of $A_+$. Then the sequence
\[
0 \rightarrow \tilde{{\cal CC}}_S(A) \stackrel {\bar{i}} {\rightarrow} 
\tilde{{\cal C}}_S(A) \stackrel {\bar{M}} {\rightarrow} 
\tilde{{\cal C}}{\cal R}_S(A) \stackrel {\bar{N}} {\rightarrow} 
\tilde{{\cal CC}}_S(A) \rightarrow 0 
\]
is exact.}

{\it The proof} is the same as that of Proposition 4 of [He2].
\hfill $\Box$

{\sc 3.4 Proposition.}~~ {\it Let $A$ be a Banach algebra with a 
bounded approximate identity $e_{\nu}, \nu \in \Lambda$,
 and let $S$ be a closed subalgebra of $A_+$. Then  ${\cal HR}^n_S(A) = \{ 0 \}$
for all $n \ge 0.$}

{\it Proof.} For $f \in  C_S^n(A)$ we define $g_{\nu} \in C^{n-1}(A)$ by
\[
 g_{\nu}(a_0,...,a_{n-1}) = f(e_\nu,a_0,...,a_{n-1}),
\]
for $a_0,...,a_{n-1} \in A$ and $\nu \in \Lambda$. 
Since the Banach space  $ C^{n-1}(A)$ is the dual space of  
$A \hat{\otimes} \dots \hat{\otimes} A$,
for the bounded net $g_{\nu}, \nu \in \Lambda$, there exists a subnet 
 $g_{\mu}, \mu \in \Lambda'$, which ${\rm weak^*}$ converges to some cochain 
 $g \in C^{n-1}(A)$. It can be checked that  $g \in C^{n-1}_S(A)$.

     For each $f \in  C_S^n(A)$ such that $\delta r^n(f) = 0$ and for 
each $\nu \in \Lambda$,
\[
0 = \delta r^n f(e_\nu,a_0,...,a_n) 
\]
\begin{equation}
\label{rdelta}
=  f(e_{\nu} a_0,a_1,...,a_n)
+ \sum_{i=0}^{n-1} (-1)^{i+1} f(e_{\nu},a_0,...,a_i a_{i+1},...,a_n).
\end{equation}
 Thus we obtain
\[
\delta r^{n-1} g(a_0,...,a_n) 
= \sum_{i=0}^{n-1} (-1)^i g(a_0,...,a_i a_{i+1},...,a_n)
\]
\[
= \lim_{\mu} \sum_{i=0}^{n-1} (-1)^i f(e_{\mu},a_0,...,a_i a_{i+1},...,a_n)
= f(a_0,a_1,...,a_n),
\]
by (\ref{rdelta}).
\hfill $\Box$

     Note that in the case where $S =  \Bbb{C} e_+$ it is easy to see that 
the statement of Proposition 3.4 is true for every Banach algebra with left 
or right bounded approximate identity. This and other conditions are given in detail in
[He2, Section 2].

{\sc 3.5 Proposition.}~~ {\it Let $A$ be a Banach algebra with a bounded approximate 
identity and let $S$ be a closed subalgebra of $A_+$. Then the $S$-relative
Connes-Tsygan exact sequence for $A$
\[
\dots  \rightarrow  {\cal H}_S^n(A)  \stackrel {B^n} {\rightarrow}  
{\cal HC}_S^{n-1}(A)  \stackrel {S^n} {\rightarrow}  {\cal HC}_S^{n+1}(A)
 \stackrel {I^{n+1}} {\rightarrow}  {\cal H}_S^{n+1}(A)   \stackrel {B^{n+1}}
 {\rightarrow} 
 {\cal HC}_S^{n}(A)  \rightarrow  \dots 
\]
exists. }

{\it Proof.} By Proposition 3.3, there are two short exact
sequences of complexes
\[
0  \rightarrow   \tilde{{\cal CC}}_S(A) \stackrel {\bar{i}} {\rightarrow} 
\tilde{{\cal C}}_S(A)  \stackrel {\bar{M}} {\rightarrow}  
\tilde{{\cal CS}}_S(A)  \rightarrow  0
\]
and 
\[
0  \rightarrow   \tilde{{\cal CS}}_S(A) \stackrel {\bar{j}} {\rightarrow} 
\tilde{{\cal CR}}_S(A)  \stackrel {\bar{N}} {\rightarrow}  
\tilde{{\cal CC}}_S(A)  \rightarrow  0
\]  
where $ \tilde{{\cal CS}}_S(A) $ is the subcomplex  ${\rm Im}(\bar{M}) =
{\rm Ker}(\bar{N}) $ of $\tilde{{\cal CR}}_S(A)$. 
Hence, by [He1, Chapter 0, Section 5.5], there exist two long exact sequences (*)
\[
\dots  \rightarrow  {\cal HC}^n_S(A) \stackrel {H^n(\bar{i})} {\rightarrow}  
{\cal H}^n_S(A)  \stackrel {H^n(\bar{M})} {\rightarrow}  {\cal HS}^n_S(A)
  \stackrel {\zeta^n} {\rightarrow} {\cal HC}^{n+1}_S(A) \stackrel {H^{n+1}(\bar{i})}
 {\rightarrow}   {\cal H}^{n+1}_S(A)  \rightarrow  \dots 
\]
and (**)
\[
\thickmuskip=.5 \thickmuskip
\medmuskip=.5 \medmuskip
\dots \rightarrow  {\cal HS}^n_S(A) \stackrel {H^n(\bar{j})}  {\rightarrow} 
{\cal HR}^n_S(A) \stackrel {H^{n+1}(\bar{N})} {\rightarrow} {\cal HC}^n_S(A) 
\stackrel {\eta^n} {\rightarrow}
  {\cal HS}^{n+1}_S(A)  \stackrel {H^{n+1}(\bar{j})} {\rightarrow}
 {\cal HR}^{n+1}_S(A)  \rightarrow  \dots. 
\]

   By Proposition 3.4, ${\cal HR}^n_S(A) = \{ 0 \}$ for all
$n \ge 0$. Thus we can see from (**) and Proposition 8 [He2] that $\eta^n$ is a 
topological isomorphism, and so that
$   {\cal HC}_S^n(A)  =  {\cal HS}_S^{n+1}(A)$
 for all $n \ge 0$.
By setting  $  {\cal HC}_S^{n-1}(A)$ instead of $ {\cal HS}_S^n(A)$  in (*), 
we get the required exact sequence, where $I^n = H^n(\bar{i}), 
S^n = \zeta^{n-1} \circ \eta^n$ and 
$B^n = (\eta^n)^{-1} \circ H^n(\bar{M})$.
\hfill $\Box$

{\sc 3.6 Proposition.}~~ {\it Let $A$ be a Banach algebra with a bounded approximate 
identity and let $S$ be a closed subalgebra of $A_+$. Then the inclusion morphism of
cochain objects $ C_S^n(A,A^*)  \rightarrow C^n(A,A^*) $ induces 
a morphism of Connes-Tsygan exact sequences for $A$, that is, a commutative diagram

\[
\arraycolsep=1pt
\begin{array}{ccccccccccccc}
\dots & \rightarrow & {\cal H}_S^n(A) & \stackrel {B^n} {\rightarrow} & 
{\cal HC}_S^{n-1}(A) & \stackrel {S^n} {\rightarrow} & {\cal HC}_S^{n+1}(A)
& \stackrel {I^{n+1}} {\rightarrow} & {\cal H}_S^{n+1}(A) &  \stackrel {B^{n+1}}
 {\rightarrow} &
 {\cal HC}_S^{n}(A) & \rightarrow & \dots 
\\
~ & ~ &  \downarrow & ~ &   \downarrow & ~ &  \downarrow & ~ & 
 \downarrow & ~ &   \downarrow & ~ & ~
\\
\dots & \rightarrow & {\cal H}^n(A) & \stackrel {B^n} {\rightarrow} & 
{\cal HC}^{n-1}(A) & \stackrel {S^n} {\rightarrow} & {\cal HC}^{n+1}(A)
& \stackrel {I^{n+1}} {\rightarrow} & {\cal H}^{n+1}(A) &  \stackrel {B^{n+1}}
 {\rightarrow} &
 {\cal HC}^{n}(A) & \rightarrow & \dots  .
\end{array}
\] }

{\it Proof.} Note that the inclusion morphism of cochain objects
$C_S^n(A, A^*) \rightarrow  C^n(A, A^*)$ gives morphisms of two pairs short exact
sequences of complexes{\hfuzz=2pt
\[
\begin{array}{ccccccccc}
0 & \rightarrow &  \tilde{{\cal CC}}_S(A) & \stackrel {\bar{i}} {\rightarrow} & 
\tilde{{\cal C}}_S(A) & \stackrel {\bar{M}} {\rightarrow} & 
\tilde{{\cal CS}}_S(A)  & \rightarrow & 0  
\\
~ & ~ &  \downarrow & ~ &   \downarrow & ~ &   \downarrow & ~ & ~
\\
0 & \rightarrow &  \tilde{{\cal CC}}(A) & \stackrel {\bar{i}} {\rightarrow} & 
\tilde{{\cal C}}(A) & \stackrel {\bar{M}} {\rightarrow} & 
\tilde{{\cal CS}}(A) & \rightarrow & 0  
\end{array}
\]}%
and
\[
\begin{array}{ccccccccc}
0 & \rightarrow &  \tilde{{\cal CS}}_S(A) & \stackrel {\bar{j}} {\rightarrow} & 
\tilde{{\cal CR}}_S(A) & \stackrel {\bar{N}} {\rightarrow} & 
\tilde{{\cal CC}}_S(A) & \rightarrow & 0  
\\
~ & ~ &  \downarrow & ~ &   \downarrow & ~ &   \downarrow & ~ & ~
\\
0 & \rightarrow &  \tilde{{\cal CS}}(A) & \stackrel {\bar{j}} {\rightarrow} & 
\tilde{{\cal CR}}(A) & \stackrel {\bar{N}} {\rightarrow} & 
\tilde{{\cal CC}}(A) &  \rightarrow & 0  .
\end{array}
\]
  
    By the cohomology analogue of Proposition II.4.2 [Ma], a morphism of two short exact
sequences of complexes
induces a morphism of long exact cohomology sequences.
Hence we have two commutative diagrams (*)
\[
\arraycolsep=1pt
\begin{array}{ccccccccccccc}
\dots & \rightarrow & {\cal HC}^n_S(A) & \rightarrow & 
{\cal H}^n_S(A) & \rightarrow & {\cal HS}^n_S(A)
& \rightarrow & {\cal HC}^{n+1}_S(A) &  \rightarrow &
 {\cal H}^{n+1}_S(A) & \rightarrow & \dots 
\\
~ & ~ &  \downarrow & ~ &   \downarrow & ~ &  \downarrow & ~ & 
 \downarrow & ~ &   \downarrow & ~ & ~
\\
\dots & \rightarrow & {\cal HC}^n(A) & \rightarrow & 
{\cal H}^n(A) & \rightarrow & {\cal HS}^n(A)
& \rightarrow & {\cal HC}^{n+1}(A) &  \rightarrow &
 {\cal H}^{n+1}(A) & \rightarrow & \dots 
\end{array}
\]
and (**)
\[
\arraycolsep=1pt
\begin{array}{cccccccccccc}
\dots & \rightarrow & {\cal HS}^n_S(A) & \rightarrow & 
{\cal HR}^n_S(A) & \rightarrow & {\cal HC}^n_S(A)
& \rightarrow & {\cal HS}^{n+1}_S(A) &  \rightarrow &
 {\cal HR}^{n+1}_S(A) &  \dots 
\\
~ & ~ &  \downarrow & ~ &   \downarrow & ~ &  \downarrow & ~ & 
 \downarrow & ~ &   \downarrow  & ~
\\
\dots & \rightarrow & {\cal HS}^n(A) & \rightarrow & 
{\cal HR}^n(A) & \rightarrow & {\cal HC}^n(A)
& \rightarrow & {\cal HS}^{n+1}(A) &  \rightarrow &
 {\cal HR}^{n+1}(A)  & \dots .
\end{array}
\]
By Proposition 3.4, ${\cal HR}^n_S(A) = {\cal HR}^n(A) = \{ 0 \}$ for all
$n \ge 0$. Thus we can see from (**) and Propositon 8 [He2] that 
there exists a commutative diagram
\[
\begin{array}{ccc}
   {\cal HC}_S^n(A) & = & {\cal HS}_S^{n+1}(A)
\\
 \downarrow & ~ &   \downarrow
\\
  {\cal HC}^n(A) & = & {\cal HS}^{n+1}(A)
\end{array}
\]
 for all $n \ge 0$.
By setting  $  {\cal HC}_S^{n-1}(A)$ instead of $ {\cal HS}_S^n(A)$ and
 $  {\cal HC}^{n-1}(A)$ instead of $ {\cal HS}^n(A)$ in (*), 
we get the required commutative diagram.
\hfill $\Box$

{\sc 3.7 Proposition.}~~ {\it Let $A$ and $D$ be Banach algebras with right or left
 bounded approximate 
identities. Suppose there exists a continuous homomorphism 
$\kappa : A \rightarrow D$. Then the associated morphism of
cochain objects $C^n(D,D^*)  \rightarrow  C^n(A,A^*) $ induces 
a morphism of Connes-Tzygan exact sequences for $A$, that is, a commutative diagram
\[
\arraycolsep=1pt
\begin{array}{cccccccccccc}
\dots & \rightarrow & {\cal H}^n(D) & \stackrel {B^n} {\rightarrow} & 
{\cal HC}^{n-1}(D) & \stackrel {S^n} {\rightarrow} & {\cal HC}^{n+1}(D)
& \stackrel {I^{n+1}} {\rightarrow} & {\cal H}^{n+1}(D) &  \stackrel {B^{n+1}}
 {\rightarrow} &
 {\cal HC}^{n}(D) &  \dots 
\\
~ & ~ &  \downarrow & ~ &   \downarrow & ~ &  \downarrow & ~ & 
 \downarrow & ~ &  \downarrow  & ~ 
\\
\dots & \rightarrow & {\cal H}^n(A) & \stackrel {B^n} {\rightarrow} & 
{\cal HC}^{n-1}(A) & \stackrel {S^n} {\rightarrow} & {\cal HC}^{n+1}(A)
& \stackrel {I^{n+1}} {\rightarrow} & {\cal H}^{n+1}(A) &  \stackrel {B^{n+1}}
 {\rightarrow} &
 {\cal HC}^{n}(A)  & \dots  .
\end{array}
\]}

{\it The proof} requires only minor modifications of that of Proposition 3.6.
\hfill $\Box$

{\sc 3.8 Proposition.}~~ {\it Let $A$ be a Banach algebra for which 
one of the following conditions is satisfied:

 {\rm (i)} $A$ has a left or right bounded approximate identity;

 {\rm (ii)} $A$ coincides with the topological square $ \overline{A^2}$ of $A$, and
that either $A$ is a flat right Banach $A$-module or ${\Bbb C}$ is a flat left Banach
$A$-module.

 Then the vanishing of ${\cal HC}^n(A)$ for 
all  $n \ge 0$ is equivalent to the vanishing of the Banach cyclic homology
 ${\cal HC}_n(A)$ for all  $n \ge 0$.}

   The definition of flat module can be found in [He1], and the definition of Banach 
cyclic
homology of a Banach algebra can be found, for example, in [He2; Section 5].

{\it Proof.} The assumption gives the existence the Connes-Tsygan 
exact sequence for cohomology of $A$ [He2, Theorems 15, 16]. We can see from this 
exact sequence that the vanishing of 
${\cal HC}^n(A)$ for all  $n \ge 0$ is equivalent to the vanishing of the 
simplicial cohomology ${\cal H}^n(A, A^*)$ for all  $n \ge 0$. 
The latter relation is equvalent to the vanishing of the simplicial homology
 ${\cal H}_n(A, A)$ for all  $n \ge 0$ [He1; Proposition 2.5.28]. 
 As was noted in [He2; Section 5], for the given assumption,  
 the canonical Connes-Tsygan exact sequence for homology of $A$ also exists. 
Thus it is easy to see from the Connes-Tsygan exact sequence for homology of $A$ 
that all these relations are  equivalent to 
the vanishing of the cyclic homology ${\cal HC}_n(A)$ for all  $n \ge 0$.
\hfill $\Box$

\section{Relative cyclic cohomology of Banach algebras}

{\sc 4.1 Theorem.}~~ {\it Let $A$ be a Banach algebra with a bounded approximate 
identity, and let $B$ be an amenable closed subalgebra of $A$. Then
\[
{\cal HC}^n(A) = {\cal HC}^n_B(A)
\] 
for all $n \ge 0$.}

{\it Proof.} Consider the commutative diagram of Proposition 3.6
\[
\arraycolsep=1pt
\begin{array}{cccccccccccccc}
0 & \rightarrow & {\cal HC}_B^0(A) & \stackrel {I^0} {\rightarrow} & 
{\cal H}_B^{0}(A) & \stackrel {B^0} {\rightarrow} & {\cal HC}_B^{-1}(A)
& \stackrel {S^{0}} {\rightarrow} & {\cal HC}_B^{1}(A) &  \stackrel {I^{1}}
 {\rightarrow} &
 {\cal H}_B^{1}(A)  & \stackrel {B^1} {\rightarrow} & {\cal HC}_B^{0}(A)
 & \dots
\\
~ & ~ &  \downarrow {\cal G}_0 & ~ &   \downarrow {\cal F}_0 & ~ &  \downarrow {\cal
G}_{-1} & ~ & 
 \downarrow  {\cal G}_1 & ~ &   \downarrow {\cal F}_1 & ~ & \downarrow {\cal G}_0  & ~ 
\\
0 & \rightarrow & {\cal HC}^0(A) & \stackrel {I^0} {\rightarrow} & 
{\cal H}^{0}(A) & \stackrel {B^0} {\rightarrow} & {\cal HC}^{-1}(A)
& \stackrel {S^{0}} {\rightarrow} & {\cal HC}^{1}(A) &  \stackrel {I^{1}}
 {\rightarrow} &
 {\cal H}^{1}(A)  & \stackrel {B^1} {\rightarrow} & {\cal HC}^{0}(A)
 & \dots 
\end{array}
\]

\[
\arraycolsep=1pt
\begin{array}{ccccccccccccc}
\dots & \rightarrow & {\cal H}_B^n(A) & \stackrel {B^n} {\rightarrow} & 
{\cal HC}_B^{n-1}(A) & \stackrel {S^n} {\rightarrow} & {\cal HC}_B^{n+1}(A)
& \stackrel {I^{n+1}} {\rightarrow} & {\cal H}_B^{n+1}(A) &  \stackrel {B^{n+1}}
 {\rightarrow} &
 {\cal HC}_B^{n}(A) & \rightarrow & \dots 
\\
~ & ~ &  \downarrow {\cal F}_n & ~ &   \downarrow {\cal G}_{n-1} & ~ & 
 \downarrow {\cal G}_{n+1} & ~ & 
 \downarrow  {\cal F}_{n+1}  & ~ &   \downarrow {\cal G}_n & ~ & ~
\\
\dots & \rightarrow & {\cal H}^n(A) & \stackrel {B^n} {\rightarrow} & 
{\cal HC}^{n-1}(A) & \stackrel {S^n} {\rightarrow} & {\cal HC}^{n+1}(A)
& \stackrel {I^{n+1}} {\rightarrow} & {\cal H}^{n+1}(A) &  \stackrel {B^{n+1}}
 {\rightarrow} &
 {\cal HC}^{n}(A) & \rightarrow & \dots  .
\end{array}
\]
Note that ${\cal HC}^0(A) = {\cal HC}_B^0(A) = A^{tr}$.
By Theorem 1.6,
\[
{\cal H}^n(A) = {\cal H}^n(A,A^*) = {\cal H}_B^n(A,A^*) = {\cal H}_B^n(A)
\] 
for all $n \ge 0$.
   
     As an induction hypothesis suppose that the vertical map 
\[
 {\cal G}_i: {\cal HC}_B^i(A) \rightarrow  {\cal HC}^i(A)
\]
is an isomorphism for each $i \le n$. Then it follows from the five lemma [Ma; Lemma
1.3.3]
that the middle vertical map ${\cal G}_{n+1}$ above is an isomorphism. Hence by Lemma
0.5.9 [He1],  ${\cal G}_{n+1}$ is topological isomorphism.
\hfill $\Box$
  
     We note that the assumptions of Theorem 4.1 are obviously satisfied by all
$C^*$-algebras $A$ and all nuclear $C^*$-subalgebras $B$ of $A_+$. Recall that all
nuclear $C^*$-algebras and only these  $C^*$-algebras are amenable. For example, all
GCR-algebras, in particular, all commutative $C^*$-algebras are amenable. Other
examples are given by the Banach algebra $A={\cal B}(E)$ and $B={\cal K}(E)$,
 where $E$ is a Banach space with the property $(\Bbb{A})$ (see [GJW]).
     
    We noted in remark after Theorem 2.4 that ${\cal K}(E)^{tr} = \{ 0 \}$ for a Banach
space  with the property $(\Bbb{A})$. Thus one can see from the following 
theorem that ${\cal HC}^n({\cal B}(E)) = 
{\cal HC}^n({\cal B}(E)/{\cal K}(E))$  for all $n \ge 0$. In particular, for 
an infinite-dimensional Hilbert space 
$H$, we have $
{\cal HC}^n({\cal B}(H)/{\cal K}(H)) = 
{\cal HC}^n({\cal B}(H)) = \{ 0 \} $ 
for all $n \ge 0$ by Theorem 4.1 [ChS1], since $ ({\cal B}(H))^{tr} = \{ 0 \}$ 
by [Hal]. Hence, by Proposition 3.8, the Banach cyclic homology of these algebras
 ${\cal HC}_n({\cal B}(H)/{\cal K}(H)) = {\cal HC}_n({\cal B}(H)) = \{ 0 \}$ for all  
$n \ge 0$.

{\sc 4.2 Theorem.}~~ {\it Let $A$ be a Banach algebra with a right or left bounded
approximate identity and let $I$ be a 
closed two-sided ideal of $A$. Suppose that $I$ is an amenable Banach algebra. 
 Then

   ${\rm (i)}$ for all even $n \ge 0$ the natural map from $ {\cal HC}^n(A/I)$ into 
${\cal HC}^n(A)$ is injective, and for all odd $n \ge 1$ the natural map from 
 $ {\cal HC}^n(A/I)$ into ${\cal HC}^n(A)$ is surjective; 

    ${\rm (ii)}$ if $I^{tr} = \{ 0 \} $ (that is, ${\rm Cen}_I I^* = \{ 0 \}$)
then
\[
{\cal HC}^n(A) = {\cal HC}^n(A/I) 
\] 
for all $n \ge 0$.} 

{\it Proof.} Consider the commutative diagram of Proposition 3.7 for the Banach algebras
$A$ and $A/I$
\[
\arraycolsep=1pt
\begin{array}{ccccccccccccc}
0 & \rightarrow & {\cal HC}^0(A/I) & \stackrel {I^0} {\rightarrow} & 
{\cal H}^{0}(A/I) & \stackrel {B^0} {\rightarrow} & {\cal HC}^{-1}(A/I)
& \stackrel {S^{0}} {\rightarrow} & {\cal HC}^{1}(A/I) &  \stackrel {I^{1}}
 {\rightarrow} &
 {\cal H}^{1}(A/I)  & \stackrel {B^1} {\rightarrow} & {\cal HC}^{0}(A/I)
\\
~ & ~ &  \downarrow {\cal G}_0 & ~ &   \downarrow {\cal L}_0 & ~ &  
\downarrow {\cal G}_{-1} & ~ & 
 \downarrow {\cal G}_1 & ~ &  \downarrow  {\cal L}_1 & ~ &  \downarrow {\cal G}_0 
\\
0 & \rightarrow & {\cal HC}^0(A) & \stackrel {I^0} {\rightarrow} & 
{\cal H}^{0}(A) & \stackrel {B^0} {\rightarrow} & {\cal HC}^{-1}(A)
& \stackrel {S^{0}} {\rightarrow} & {\cal HC}^{1}(A) &  \stackrel {I^{1}}
 {\rightarrow} &
 {\cal H}^{1}(A)  & \stackrel {B^1} {\rightarrow} & {\cal HC}^{0}(A) 
\end{array}
\]
\[
\arraycolsep=1pt
\begin{array}{cccccccccccc}
\dots & \rightarrow & {\cal H}^n(A/I) & \stackrel {B^n} {\rightarrow} & 
{\cal HC}^{n-1}(A/I) & \stackrel {S^n} {\rightarrow} & {\cal HC}^{n+1}(A/I)
& \stackrel {I^{n+1}} {\rightarrow} & {\cal H}^{n+1}(A/I) &  \stackrel {B^{n+1}}
 {\rightarrow} &
 {\cal HC}^{n}(A/I) & \dots 
\\
~ & ~ &  \downarrow  {\cal L}_n & ~ &   \downarrow  {\cal G}_{n-1} & ~ &  
\downarrow  {\cal G}_{n+1} & ~ & 
 \downarrow  {\cal L}_{n+1} & ~ &  \downarrow  {\cal G}_n & ~  
\\
\dots & \rightarrow & {\cal H}^n(A) & \stackrel {B^n} {\rightarrow} & 
{\cal HC}^{n-1}(A) & \stackrel {S^n} {\rightarrow} & {\cal HC}^{n+1}(A)
& \stackrel {I^{n+1}} {\rightarrow} & {\cal H}^{n+1}(A) &  \stackrel {B^{n+1}}
 {\rightarrow} &
 {\cal HC}^{n}(A) & \dots  .
\end{array}
\]

    Note that the maps ${\cal G}_0: {\cal HC}^0(A/I)
 \rightarrow  {\cal HC}^0(A): f  \mapsto f \circ g $ and $ {\cal L}_0: {\cal H}^0(A/I)
 \rightarrow  {\cal H}^0(A): f  \mapsto f \circ g $ are injective. By Theorem 2.4,
\[
{\cal H}^n(A) = {\cal H}^n(A,A^*) = {\cal H}^n(A/I,(A/I)^*) = {\cal H}^n(A/I),
\]
so that $ {\cal L}_n$ is an isomorphism, 
for all $n \ge 2$ and the natural  ${\cal L}_1$ map from 
$$ {\cal H}^1(A/I,(A/I)^*)$$ into 
${\cal H}^1(A,A^*)$ is surjective. Then it follows
 from the five lemma
that the middle vertical map ${\cal G}_1$ above is surjective.

     Suppose, inductively, that the vertical map, for each $i \le n$, 
\[
 {\cal G}_i: {\cal HC}^i(A/I) \rightarrow  {\cal HC}^i(A)
\]
is injective if $i$ is even and surjective if $i$ is odd. The result then follows 
from the five lemma.
\hfill $\Box$

  Note that in particular Theorem 4.2 applies whenever $A$ is a $C^*$-algebra and $I$ is
an amenable  closed ideal.

{\sc 4.3 Proposition.}~~ {\it Let $A$ and $D$ be Banach algebras with right or left
 bounded approximate 
identities, and let $\kappa : A \rightarrow D$ be a continuous homomorphism.
If $\kappa$ induces a topological isomorphism   
\[
{\cal H}^n(D, D^*) \rightarrow {\cal H}^n(A, A^*) 
\] 
 for all $n \ge 0$, then it  induces a topological isomorphism 
\[
{\cal HC}^n(D)  \rightarrow {\cal HC}^n(A) 
\] 
for all $n \ge 0$, and conversely.}

{\it Proof.} The forward implication is a repetition of that of Theorem 4.1 with the
commutative diagram of Proposition 3.7. The converse statement follows easily from the
five lemma.
\hfill $\Box$

{\sc 4.4 Proposition.}~~ {\it Let $A_i$ be a Banach algebra with 
identity $e_i, i = 1, ..., m,$ and let $A$ be the Banach algebra direct sum 
$ \bigoplus_{i=1}^m A_i$.
Then  
\[ 
{\cal HC}^n(A) =  \bigoplus_{i=1}^m {\cal HC}^n(A_i) 
\]
for all $n \ge 0$.}

{\it Proof.} By Theorem 4.1,   ${\cal HC}^n(A) = {\cal HC}_B^n(A)$, 
where $n \ge 0$ and $B$ is the Banach subalgebra of $A$ generated by 
$ \{e_i, i = 1, ..., m\}$.
By Proposition 1.7, the canonical projections from $A$ to $A_i, i = 1, ..., m,$ 
induce a topological isomorphism of complexes ${\cal C}_B(A,A^*)  \rightarrow  
\bigoplus_{i=1}^m {\cal C}(A_i, A_i^*)$. Thus 
${\cal HC}_B^n(A) =  \bigoplus_{i=1}^m {\cal HC}^n(A_i)$.
The result now follows directly. 
\hfill $\Box$

{\sc 4.5 Corollary.}~~ {\it Let ${\cal R}$ be a von Neumann algebra, let
\[
{\cal R} = {\cal R}_{I_f}  \oplus {\cal R}_{I_\infty}  
\oplus {\cal R}_{II_1} \oplus {\cal R}_{II_\infty} 
\oplus {\cal R}_{III}
\]
be the central direct summand decomposition of ${\cal R}$ into von Neumann algebras of 
types 
$I_f, I_{\infty}, II_1, II_{\infty}, III$.
Then  
\[ 
{\cal HC}^n({\cal R}) =
 {\cal HC}^n({\cal R}_{I_f})  \oplus
 {\cal HC}^n({\cal R}_{II_1})  
\]
for all $n \ge 0$.}

{\it Proof.} As we noted in the proof of Proposition 1.8, there are no non-zero bounded
traces on $ {\cal R}_{I_\infty}, {\cal R}_{II_\infty}$ and $ {\cal R}_{III}$. Thus, by
 Theorem 4.1 [ChS1], their Banach cyclic cohomology groups vanish for all  $n \ge 0$.
\hfill $\Box$

   Note that, by Theorem 25 [He2], for a von Neumann algebra ${\cal R}_{I_m}$ of type
 $I_m$, where $m < \infty$, we have  ${\cal HC}^n({\cal R}_{I_m}) = 
{\cal R}_{I_m}^{tr}$ for all even $n$, and   ${\cal HC}^n({\cal R}_{I_m}) = \{ 0 \}$ for all 
odd $n$.

{\sc 4.6 Corollary.}~~ {\it Let $A$ be a $C^*$-algebra without non-zero bounded traces.
Then the Banach cyclic homology
 ${\cal HC}_n(A) = \{ 0 \}$ for all  $n \ge 0$.}

{\it Proof.} It follows from Proposition 3.8 and Theorem 4.1 [ChS1].
\hfill $\Box$ 
 
  Note that in particular Corollary 4.6 applies whenever $A$ is a properly infinite von
Neumann algebra (see Proposition 1.8), or a stable $C^*$-algebra [Fa; Theorem 1.1].
 The calculation of the (Banach) cyclic homology groups for stable
$C^*$-algebras, for ${\cal B}(H)$ and for the Calkin algebra on a 
separable Hilbert space was given in [Wo].

     Now we give the Banach version of the algebraic statement of Theorem 9 [Ka1] of
L. Kadison.

{\sc 4.7 Proposition.}~~ {\it Let $A_1$ and $A_2$ be unital Banach algebras, 
let $Y$ be a unital Banach $A_1-A_2$-bimodule, and let   
\[
{\cal U} = \left[ \begin{array}{cc}
A_1 & Y
\\
0 & A_2
\end{array} \right]
\]
the natural triangular matrix algebra.
 Then the two canonical projections from ${\cal U}$ to $A_1$ and $A_2$ induce a 
topological isomorphism
\[
{\cal HC}^n({\cal U}) = {\cal HC}^n(A_1)  \oplus {\cal HC}^n(A_2)
\] 
for all $n \ge 0$.}

{\it Proof.} By Theorem 4.1 and Proposition 1.10, we obtain
\[
{\cal HC}^n({\cal U}) = {\cal HC}_B^n({\cal U}) = 
{\cal HC}^n(A_1)  \oplus {\cal HC}^n(A_2),
\]
where the amenable subalgebra $B$ of ${\cal U}$ was defined before Proposition 1.10.
\hfill $\Box$

\end{document}